\documentclass[submission]{FPSAC2018}

\usepackage[latin1]{inputenc}
\usepackage{subfigure, amsmath, verbatim}

\usepackage{tikz}
\usetikzlibrary{shapes.geometric,positioning}
\usepackage{xcolor}
\usetikzlibrary{arrows, shapes}
\newtheorem{theorem}{Theorem}%[section]
\newtheorem{proposition}[theorem]{Proposition}

\newtheorem{corollary}[theorem]{Corollary}

\theoremstyle{definition}

\newtheorem{definition}[theorem]{Definition}

\newtheorem{example}[theorem]{Example}

\newcommand{\spe}[1]{\mathtt{#1}}

\DeclareMathOperator{\inc}{inc}

\DeclareMathOperator{\Des}{Des}

\begin{document}

\title{The Hopf monoid of Megagreedoids}

\author[J.White]{Jacob A.White\thanks{\href{mailto:jacob.white@utrgv.edu}{jacob.white@utrgv.edu}}\addressmark{1}}
\address{\addressmark{1}School of Mathematical and Statistical Sciences,
University of Texas - Rio Grande Valley,
Edinburg, TX 78539}
%\email{jwhite@math.tamu.edu}
%\urladdr{http://www.math.tamu.edu/~jwhite/}

%\subjclass[2010]{Primary 18D10, 18D35}
\received{\today}

%18D10 Monoidal categories
%18D35 Structured objects in a category

\abstract{We introduce megagreedoids, which generalize polymatroids, megamatroids, and greedoids. We define a quasisymmetric function invariant for a megagreedoid, and show that it has a positive expansion in the basis of fundamental quasisymmetric functions. Our proof involves lexicographic shellability. We also show that megagreedoids form a Hopf monoid.
A running example is a megagreedoid associated to a rooted connected graph, and the resulting generalization of the chromatic symmetric function.}

\keywords{Hopf Algebras, Combinatorial Species, Quasisymmetric Functions, Shellability, Greedoids } 

\maketitle

\section{Introduction}
In \cite{aguiar-bergeron-sottile}, Aguiar, Bergeron, and Sottile show that quasisymmetric functions were the terminal combinatorial Hopf algebra, which 'explains the ubiquity of quasi-symmetric functions as generating functions in combinatorics'. Many combinatorial Hopf algebras have been studied, such as graphs \cite{schmitt}, posets \cite{joni-rota}, building sets \cite{grujic}, simplicial complexes \cite{benedetti}, matroids \cite{billera-jia-reiner}, megamatroids \cite{derksen-fink}, and hypergraphs \cite{grujic}. Aguiar and Mahajan initiated the study of Hopf monoids in species \cite{aguiar-mahajan-1}. In many cases, the Hopf monoid/algebra has a nice basis for which the resulting quasisymmetric functions are positive when expanded in the basis of fundamental quasisymmetric functions. However, there is no general proof of this fact. Moreover, for many of the newer combinatorial Hopf algebras, $F$-positivity is unknown. 

In recent work, Aguiar and Ardila \cite{aguiar-ardila} studied the Hopf monoid of generalized permutahedra, which generalizes all of these other examples. They proved a grouping-free, cancellation-free formula for the antipode, and proved a combinatorial reciprocity result, generalizing many known results. They also study a related polynomial invariant for generalized permutahedra. Their work is highly geometric, highlighting the fact that there is quite a bit of unexplored intersection between geometric combinatorics and Hopf algebras. This paper explores this intersection further.

The goal of this paper is to introduce megagreedoids, defined in Section \ref{sec:setting}. These generalize polymatroids, megamatroids \cite{derksen-fink}, and greedoids \cite{greedoid}, which in turn generalize numerous combinatorial objects. We also study a new class of megagreedoids which come from connected rooted graphs, to give evidence that there are other interesting classes of megagreedoids that have yet to be studied.

In section \ref{sec:quasi}, we show that there is a natural quasisymmetric function invariant associated to megagreedoids, involving enumerating certain functions. Moreover, this invariant is $F$-positive, generalizing known results for the chromatic symmetric function \cite{stanley-coloring-1}, the $P$-partition generating function of Gessel \cite{gessel}, and the matroid invariant of Billera, Jia, and Reiner \cite{billera-jia-reiner}. We also obtain $F$-positivity results for new quasisymmetric functions associated to polymatroids, greedoids, and rooted connected graphs. We also state a combinatorial reciprocity result for the resulting polynomial invariant obtained through specialization.

Our approach to $F$-positivity is to define a relative simplicial complex which generalizes the coloring complex of a graph \cite{steingrimsson}. We detail this approach in Section \ref{sec:shellability}. It involves a simple extension of lexicographic shellability to relative order complexes. We use the greedy algorithm to show that the resulting relative simplicial complexes coming from megagreedoids are shellable. The proof is omitted, although it is elementary. However, the simplicity of the proofs gives evidence that the notion of megagreedoid is useful.
Thus, we see a rich interplay between topological combinatorics, optimization theory, and Hopf algebras.

In Section \ref{sec:hopf}, we review the definition of Hopf monoids and chromatic quasisymmetric functions. We show that megagreedoids form the basis of a Hopf monoid $\spe{MG}$ which contains the Hopf monoid of graphs, matroids, posets, and polymatroids as Hopf submonoids. Thus we obtain $F$-positivity for previously studied invariants. 

Admittedly, there is much left to discover, both about megagreedoids, and the chromatic quasisymmetric function of rooted connected graphs. One other note: the material in Section \ref{sec:shellability} and Section \ref{sec:hopf} can be read independently of each other, and are likely of interest to different readers. However, the former is how we prove $F$-positivity, and the latter is how we prove our combinatorial reciprocity result, so each point of view brings something different.
\section{Megagreedoids}
\label{sec:setting}

\begin{definition}
A megagreedoid on $I$ is a pair $(A,r)$, where $A \subset 2^I$, $r: A \to \mathbb{R}$, subject to:
\begin{enumerate}
\item for every $S \in A$, there exists $x\in S$, and $y \in I \setminus S$, such that $S - \{x\}, S \cup \{y\} \in A$,
\item for all $X \subset Y \in A$, and $z \in I \setminus Y$, if $X \cup \{z\}, Y \cup \{z\} \in A$, then $r(Y \cup \{z\}) - r(Y) \leq r(X \cup \{z\}) - r(X)$.
\end{enumerate}
\end{definition}

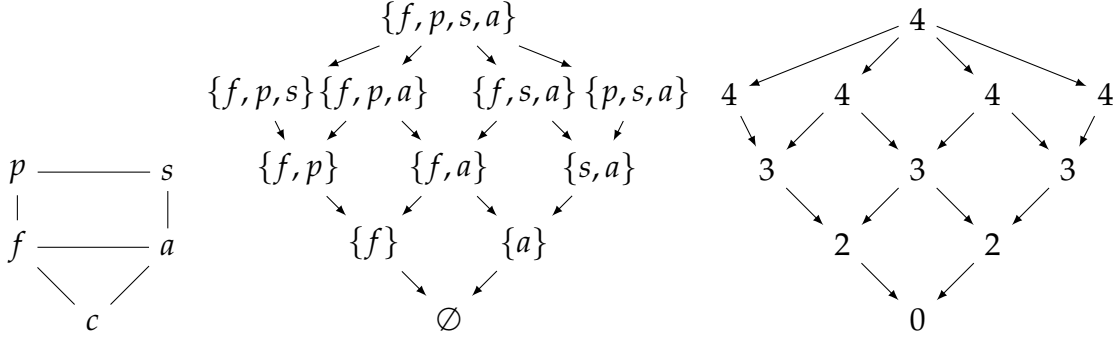
\begin{figure}[t]

\begin{center}
\begin{tikzpicture}
\node (e) at (0,0) {$c$};
\node (f) at (-1,1) {$f$};
\node (n) at (1,1) {$a$};
\node (u) at (-1,2) {$p$};
\node (s) at (1,2) {$s$};

\draw[-] (f) -- (e);
\draw[-] (n) -- (e);
\draw[-] (u) -- (f);
\draw[-] (f) -- (n);
\draw[-] (u) -- (s);
\draw[-] (s) -- (n);

\end{tikzpicture}
\begin{tikzpicture}
\node (e) at (0,0) {$\emptyset$};
\node (a) at (-1,1) {$\{f\}$};
\node (d) at (1,1) {$\{a\}$};
\node (ab) at (-2,2) {$\{f,p \}$};
\node (ad) at (0,2) {$\{f,a \}$};
\node (cd) at (2,2) {$\{s,a\}$};
\node(nd) at (-2.5,3) {$\{f,p,s \}$};
\node (nc) at (-1,3) {$\{f,p,a\}$};
\node (nb) at (1,3) {$\{f,s,a \}$};
\node (na) at (2.5,3) {$\{p,s,a \}$};
\node (t) at (0,4) {$\{f,p,s,a \}$};
\draw[->, >=latex] (a) -- (e);
\draw[->, >=latex] (d) -- (e);
\draw[->, >=latex] (ab) -- (a);
\draw[->, >=latex] (ad) -- (a);
\draw[->, >=latex] (ad) -- (d);
\draw[->, >=latex] (cd) -- (d);
\draw[->, >=latex] (nd) -- (ab);
\draw[->, >=latex] (nc) -- (ab);
\draw[->, >=latex] (nc) -- (ad);
\draw[->, >=latex] (nb) -- (ad);
\draw[->, >=latex] (nb) -- (cd);
\draw[->, >=latex] (na) -- (cd);
\draw[->, >=latex] (t) -- (na);
\draw[->, >=latex] (t) -- (nb);
\draw[->, >=latex] (t) -- (nc);
\draw[->, >=latex] (t) -- (nd);

\end{tikzpicture}
\begin{tikzpicture}
\node (e) at (0,0) {$0$};
\node (a) at (-1,1) {$2$};
\node (d) at (1,1) {$2$};
\node (ab) at (-2,2) {$3$};
\node (ad) at (0,2) {$3$};
\node (cd) at (2,2) {$3$};
\node(nd) at (-2.5,3) {$4$};
\node (nc) at (-1,3) {$4$};
\node (nb) at (1,3) {$4$};
\node (na) at (2.5,3) {$4$};
\node (t) at (0,4) {$4$};
\draw[->, >=latex] (a) -- (e);
\draw[->, >=latex] (d) -- (e);
\draw[->, >=latex] (ab) -- (a);
\draw[->, >=latex] (ad) -- (a);
\draw[->, >=latex] (ad) -- (d);
\draw[->, >=latex] (cd) -- (d);
\draw[->, >=latex] (nd) -- (ab);
\draw[->, >=latex] (nc) -- (ab);
\draw[->, >=latex] (nc) -- (ad);
\draw[->, >=latex] (nb) -- (ad);
\draw[->, >=latex] (nb) -- (cd);
\draw[->, >=latex] (na) -- (cd);
\draw[->, >=latex] (t) -- (na);
\draw[->, >=latex] (t) -- (nb);
\draw[->, >=latex] (t) -- (nc);
\draw[->, >=latex] (t) -- (nd);

\end{tikzpicture}

\end{center}
\caption{A rooted graph, its megagreedoid, and the rank function.}
\label{fig:graph}
\end{figure}

Let $G$ be a connected graph with vertex set $I \cup \{r\}$, where $r$ is a distinguished root vertex. For a subset $S \subseteq I$, let $G|_S$ denote the induced subgraph on $S$. Define $M_G = (A_G, \inc_G)$, where $A_G = \{S \subseteq I: G|_{I \cup \{r\}} \mbox{is connected} \}$. Moreover, $\inc_G(S)$ is the number of edges of $G$ with at least one endpoint in $S$, and the other in $I$. Then $M_G$ is a megagreedoid. An example is given in \ref{fig:graph}.

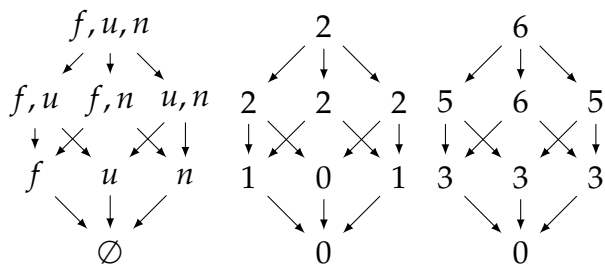
\begin{figure}[t]
\begin{center}

\begin{tikzpicture}
\node (e) at (0,0) {$\emptyset$};
\node (f) at (-1,1) {$f$};
\node (n) at (1,1) {$n$};
\node (fn) at (0,2) {$f,n $};
\node (fu) at (-1,2) {$f,u$};
\node (un) at (1,2) {$u,n $};
\node (t) at (0,3) {$f,u,n $};
\node (u) at (0,1) {$u$};

\draw[->, >=latex] (f) -- (e);
\draw[->, >=latex] (n) -- (e);
\draw[->, >=latex] (fn) -- (f);
\draw[->, >=latex] (fn) -- (n);

\draw[->, >=latex] (u) -- (e);
\draw[->, >=latex] (fu) -- (u);
\draw[->, >=latex] (un) -- (u);

\draw[->, >=latex] (t) -- (un);
\draw[->, >=latex] (t) -- (fu);

\draw[->, >=latex] (fu) -- (f);
\draw[->, >=latex] (un) -- (n);
\draw[->, >=latex]  (t) -- (fn);
\end{tikzpicture}
\begin{tikzpicture}
\node (e) at (0,0) {$0$};
\node (f) at (-1,1) {$1$};
\node (n) at (1,1) {$1$};
\node (fn) at (0,2) {$2 $};
\node (fu) at (-1,2) {$2$};
\node (un) at (1,2) {$2 $};
\node (t) at (0,3) {$2 $};
\node (u) at (0,1) {$0$};

\draw[->, >=latex] (f) -- (e);
\draw[->, >=latex] (n) -- (e);
\draw[->, >=latex] (fn) -- (f);
\draw[->, >=latex] (fn) -- (n);

\draw[->, >=latex] (u) -- (e);
\draw[->, >=latex] (fu) -- (u);
\draw[->, >=latex] (un) -- (u);

\draw[->, >=latex] (t) -- (un);
\draw[->, >=latex] (t) -- (fu);

\draw[->, >=latex] (fu) -- (f);
\draw[->, >=latex] (un) -- (n);
\draw[->, >=latex]  (t) -- (fn);
\end{tikzpicture}
\begin{tikzpicture}
\node (e) at (0,0) {$0$};
\node (f) at (-1,1) {$3$};
\node (n) at (1,1) {$3$};
\node (fn) at (0,2) {$6 $};
\node (fu) at (-1,2) {$5$};
\node (un) at (1,2) {$5 $};
\node (t) at (0,3) {$6 $};
\node (u) at (0,1) {$3$};

\draw[->, >=latex] (f) -- (e);
\draw[->, >=latex] (n) -- (e);
\draw[->, >=latex] (fn) -- (f);
\draw[->, >=latex] (fn) -- (n);

\draw[->, >=latex] (u) -- (e);
\draw[->, >=latex] (fu) -- (u);
\draw[->, >=latex] (un) -- (u);

\draw[->, >=latex] (t) -- (un);
\draw[->, >=latex] (t) -- (fu);

\draw[->, >=latex] (fu) -- (f);
\draw[->, >=latex] (un) -- (n);
\draw[->, >=latex]  (t) -- (fn);
\end{tikzpicture}

\end{center}
\caption{A boolean lattice, a greedoid, and a polymatroid.}
\label{fig:greed}
\end{figure}

Megagreedoids generalize polymatroids, which are megagreedoids where $A = 2^I$. These are also called generalized permutahedra. An example is given in Figure \ref{fig:greed}. When $A$ is the set of lower order ideals of a poset, then $M$ is equivalent to a megamatroid \cite{derksen-fink}.

Megagreedoids also generalize greedoids \cite{greedoid}.
A \emph{greedoid} is a function $r: I \to \mathbb{N}$, subject to:
\begin{enumerate}
\item $r(A) \leq |A|$ for all $A \subseteq I$,
\item $r(A) \leq r(B)$ for all $A \subseteq B$,
\item if $r(A) = r(A \cup \{x\}) = r(A \cup \{y\})$ for $x, y \in I \setminus A$, then $r(A) = r(A \cup \{x,y\})$.
\end{enumerate}
For greedoids, the greedy algorithm minimizes \emph{rank-feasible} linear functions \cite{greedoid}.
A subset $S \subseteq I$ is \emph{rank feasible} if $r(S \cup X) \leq r(S) + |X|$ for all $X \subset I \setminus S$. Given a greedoid $G$ with rank function $r$, let $M_G = (A_G, r)$, where $A_G$ is the collection of rank-feasible subsets of $G$. Then $M_G$ is a megagreedoid. In Figure \ref{fig:greed}, the middle graph is a greedoid: the megagreedoid is obtained by deleting $\{u\}$, which has rank $0$.

An important class of greedoids comes from posets. Given a poset $P$ on $I$, and $S \subseteq I$, let $r(S)$ be the maximum size of an order ideal of $P$ contained in $S$. Then $r: 2^I \to \mathbb{N}$ is a greedoid.
Equivalently, the megagreedoid of $P$, $M_P = (A_P, r)$, is given by setting $A_P$ equal to the collection of lower order ideals of $P$, and letting $r$ be the usual cardinality function.

\begin{comment}
A new class of megagreedoid come from \emph{posets with inconsistent pairs}, which is a combinatorial structure defined by Ardila, Owen, and Sullivant in their study of CAT(0) spaces. Recall that a poset with inconsistent pairs is a poset $P$ with a set $I$ of edges on $V$, subject to:
\begin{enumerate}
\item If $\{p, q\} \in I$, then there does not exist $r \in P$ with $p \leq r$, $q \leq r$.
\item If $\{p, q \} \in I$, and $p \leq p'$, $q \leq q'$, then $\{p', q' \} \in I$.
\end{enumerate}
Given a poset with inconsistent pairs, let $A_P$ be the collection of order ideals. Given an order ideal $S$, let $\inc(S)$ denote the number of inconsistent pairs with at least one element in $S$. Then $(A_P, \inc)$ is a megagreedoid.

As one final example of a megagreedoid, let $n \in \mathbb{N}$, and let $A_n = \{ Inv(\sigma): \sigma \in \mathfrak{S}_n \}$, the collection of inversion sets of permutations. Let $r(S) = |S|$. Then $(A_n, r)$ is a megagreedoid. $A_n$ is also known as the right weak order on the symmetric group. This construction holds for aribtrary finite Coxeter groups, and deserves further study. At the very least, it represents another example of a megagreedoid that is tangential to the others considered in this paper.
\end{comment}
We define the direct sum of two megagreedoids. Given $I = S \sqcup T$, and megagreedoids $(A,r)$ on $S$, $(B,s)$ on $T$, let $(A,r) \cdot (B,s) = (A \cdot B, r \cdot s)$, where $A \cdot B = \{X \cup Y: X \in A, Y \in B \}$, and $r \cdot s: A \cdot B \to \mathbb{R}$ is given by $(r \cdot s)(X) = r(X \cap S) + s(X \cap T)$.

Given $I$, and a megagreedoid $M=(A,r)$ on $I$, let $S \in A$. We define the \emph{restriction} $M|_S = (A|_S, r|_S)$, where $A|_S = \{X \in A: X \subseteq S \}$, and $r|_S$ is just the restriction of $r$ to $A|_S$. The \emph{contraction} $M/S$ is defined to be $M/S = (A/S, r/S)$, where $A / S = \{X \subseteq I \setminus S: X \cup S \in A \}$, and $r/S: A/S \to \mathbb{R}$ is given by $r/S(X) = r(X \cup S) - r(S)$ for all $X \in A/S$.

Define the \emph{base polyhedron} $P_M$ of $M$ to be given by inequalities $\sum_{s \in S} x_s \leq r(S)$ for all $S \in A$. When $M$ is a polymatroid, 
this is a generalized permutahedron. A function $f: I \to \mathbb{R}$ is $r$-feasible if $\{i: f(i) \leq c \} \in A$ for all $c \in \mathbb{R}$. It is strongly $r$-feasible if any sufficiently small perturbation of $f$ results in another $r$-feasible function. For example, consider the function $F: \{f,p,s,a \} \to \mathbb{R}$ given by $F(a) = F(s) = 1$, and $F(p) = F(f) = 1$. Then $F$ is $r$-feasible for the megagreedoid in Figure \ref{fig:graph}. However, it is not strongly feasible, because if the value of $F(s)$ is changed to $.9$, then the function is no longer $r$-feasible.

\section{Quasisymmetric Function Invariants}
\label{sec:quasi}

To a megagreedoid, we associate a quasisymmetric function.
Recall that a quasisymmetric function is a power series $f$ in variables $\{x_i: i \in \mathbb{N} \}$, such that, for any $i_1 < i_2 < \cdots < i_k$, and any sequence $(\alpha_1, \cdots, \alpha_k)$, the coefficient of $x_{i_1}^{\alpha_1}\cdots x_{i_k}^{\alpha_k}$ in $f$ is equal to the coefficient of $x_1^{\alpha_1} \cdots x_k^{\alpha_k}$. 

\begin{definition}
Given a megagreedoid $M = (A,r)$ on $I$, with base polyhedron $P_M$, $f:I \to \mathbb{N}$ is $r$-generic if it strongly $r$-feasible, and is minimized at a unique vertex of $P_M$.
Weighting $r$-generic functions by $\textbf{x}_f = \prod_{i \in I} x_{f(i)}$, the \emph{$r$-generic quasisymmetric function} $\chi(M)$ is:
$$\chi(M) = \sum_{f: I \to \mathbb{N}} \textbf{x}_f.$$

\end{definition}
There is also a corresponding $r$-generic polynomial, defined by $\chi(M, n) = $ the number of $r$-generic functions $f:I \to [n]$. 
In the case of matroids, $\chi(M)$ was studied by Billera, Jia, and Reiner \cite{billera-jia-reiner}. For polymatroids, $\chi(M,n)$ was studied by Aguiar and Ardila \cite{aguiar-ardila}.

Let $G$ be a rooted connected graph, with associated megagreedoid $M_G$.
Then $f:I \to \mathbb{N}$ is strongly feasible if $f^{-1}([i])$ induces a connected subgraph of $G$ for all $i \in \mathbb{N}$. Moreover, $f$ is $r$-generic if $f^{-1}(i)$ is an independent set of $G$ for all $i$. 
Thus $f$ is a proper coloring such that, for every $v \in V$, there is a path $r, v_1, \ldots, v_k = v$ with $f(v_1) < f(v_2) < \cdots < f(v_k)$. When $r$ is connected to every vertex, we obtain the chromatic symmetric function \cite{stanley-coloring-1}. Ordering the vertices as $f < p < s< a$, then for the graph in Figure \ref{fig:graph}, we see that the function $(1,1,2,3)$ is not proper, because there is an edge $fp$. On the other hand, $(2,1,1,2)$ is not proper, because the vertex $u$ is not adjacent to $r$, so it cannot be colored $1$. However, the colorings $(1,2,3,4)$ and $(1,2,3,2)$ are proper.

Now let $(P,I)$ be a poset. In this case, a function $f$ is feasible if and only if $f^{-1}([i])$ is an order ideal for all $i$. Equivalently, $f(p) \leq f(q)$ for all $p \leq q$. The requirement that $f$ be strongly feasible forces $f(p) < f(q)$ for all $p < q$.  Thus, we obtain Gessel's $P$-partition generating function.

\begin{proposition}
Let $M = (A,r)$ be a megagreedoid on $I$, and let $N$ be a megagreedoid on $J$.
\begin{enumerate}
\item $\chi( M \cdot N) = \chi(M) \chi(N).$
and 
\item $\chi(M, n+m) = \sum_{S \in A} \chi(M|_S, n) \cdot \chi(M/S, m)$
\end{enumerate}
\end{proposition}
Recall the basis of fundamental quasisymmetric functions $F_{S}, n$, defined by:

\[F_{S,n} = \sum_{\substack{i_1 \leq i_2 \leq \cdots \leq i_n \\ j \in S \to i_j < i_{j+1}}} x_{i_1} \cdots x_{i_n}.\]
Similarly, the monomial quasisymmetric functions are defined by: $F_{S,n} = \sum\limits_{T \supset S} M_{T,n}$.
Fix a linear order on $I$. Given a permutation $\sigma$ of $I$, $\sigma$ is $A$-feasible if $\{\sigma_1, \ldots, \sigma_i \} \in A$ for all $i$. Let $\mathfrak{S}(A)$ be the set of $A$-feasible permutations.
Given $\sigma \in \mathfrak{S}(A)$, there is an $r$-descent at $i$ if any of the following conditions are met:
\begin{enumerate}
\item $r(\{\sigma_1, \ldots, \sigma_i \}) > r(\{\sigma_1, \ldots, \sigma_{i+1} \})$,
\item $r(\{\sigma_1, \ldots, \sigma_i \}) = r(\{\sigma_1, \ldots, \sigma_{i+1} \})$ and $\sigma_i > \sigma_{i+1}$,
\item $\{\sigma_1, \ldots, \sigma_{i-1}, \sigma_{i+1} \} \not\in A$,
\item $r(\{\sigma_1, \ldots, \sigma_i \}) - r(\{\sigma_1, \ldots, \sigma_{i-1} \}) > r(\{\sigma_1, \ldots, \sigma_{i+1} \}) - r(\{\sigma_1, \ldots, \sigma_{i-1}, \sigma_{i+1} \})$
\end{enumerate}
We let $\Des(\sigma)$ be the set of $r$-descents.
\begin{theorem}
Let $M = (A,r)$ be a megagreedoid on $I$. Then
\[\chi(A,r) = \sum_{\sigma \in \mathfrak{S}(A)} F_{\Des(\sigma), |I|}.\] 
\end{theorem}
Note that, every graph $G$ can be viewed as a rooted connected graph $G'$, by adding a new root vertex $r$ that is attached to all other vertices. Then $\chi(M_{G'}) = \chi(G)$, so we obtain $F$-positivity of the chromatic symmetric function of a graph. Our notion of descent is different from the one used by Steingr\'imsson \cite{steingrimsson}. Also, since every matroid is also a megagreedoid, we get $F$-positivity of the Billera-Jia-Reiner invariant \cite{billera-jia-reiner}.

Now we do some example calculations. Let $G$ be the rooted graph from Figure \ref{fig:graph}, and let $M_G = (A,r)$ be the associated megagreedoid. We order $\{f,p,s,a\}$ by $a < f < p < s$. Then $asfp \in \mathfrak{S}(A)$ and $\Des(asfp) = \{1,2,3\}$. 
Then $\chi(M_G) = 6F_{\{1,2,3\},4}+2F_{\{1,3\},4}$.
For the greedoid $G$ from Figure \ref{fig:greed}, $\chi(M_G) = F_{\Des(fun), 3} + F_{\Des(fnu), 3} + F_{\Des(nfu), 3} + F_{\Des(nuf), 3} = F_{\{1,2\},3} + F_{\{2\},3}+F_{\{1,2\},3}+ F_{\{1,2\},3} = F_{\{2\},3} + 3F_{\{1,2\},3}$.
For the polymatroid $M$ in \ref{fig:greed}, $\chi(M) = F_{\{2\}, 3} + 5F_{\{1,2\},3}$.

We also obtain the following combinatorial reciprocity result:
\begin{theorem}
Let $M$ be a megagreedoid. Then
$(-1)^{|I|}\chi(M, -n) = \sum_{f:I \to [n]} v(f)$
where $f$ is $r$-feasible, and $v(f)$ is the number of vertices of $P_M$ where $f$ is minimized.
\end{theorem}

\begin{corollary}
Let $G$ be a connected rooted graph, and let $M_G$ be the associated megagreedoid. Then
$(-1)^{|I|}\chi(M_G, -1)$ is the number of acyclic orientations of $G$ that contain a directed spanning tree rooted at $r$, which is a sink.
\end{corollary}

\section{Shellability of Relative Complexes}
\label{sec:shellability}
In this section, we prove $F$-positivity for $\chi(M)$ by first relating our invariant to a face enumerator for a balanced relative simplicial complex $(\Sigma(M), \Gamma(M))$.  Then we discuss shellability of relative simplicial complexes. Then we show that $(\Sigma(M), \Gamma(M))$ is shellable, by using a greedy algorithm. This is equivalent to a generalization of lexicographic shellability, although our approach requires less terminology. In the full version, we will discuss the generalization of lexicographic shellability, as it may be of independent interest.

A \emph{relative simplicial complex} on a finite set $X$ is a pair $\Gamma \subset \Sigma$ of simplicial complexes. Given a relative complex $\Gamma \subset \Sigma$, a cocomplex $\Psi$ is the set of faces $\Delta \setminus \Gamma$.
A \emph{balanced} relative simplicial complex of dimension $d$ is a triple $(\Sigma, \Gamma, \rho)$ where $\rho: V \to [d]$, such that $\rho$ is a proper coloring of the one-skeleton of $\Sigma$. The associated quasisymmetric function is given by 
\[F(\Sigma, \Gamma, \rho) = \sum_{\sigma \in \Sigma \setminus \Gamma} M_{\rho(\sigma), d} \]
where $\rho(\sigma) = \{\rho(x): x \in \sigma \}$.
Given a megagreedoid $M = (A,r)$, let $\Sigma(M)$ be the order complex of $A$, ordered by inclusion. Let $\Gamma_r(M)$ be defined so that $\Sigma(M) \setminus \Gamma_r(M)$ are those chains $S_1 \subset S_2 \subset \cdots \subset S_k$ where, for every $i$, $[S_i, S_{i+1}]$ is a boolean interval, and $r|_{S_{i+1}} / S_i$ is modular. Also, we define $\ell(C)$ of a chain to be the length of the chain.

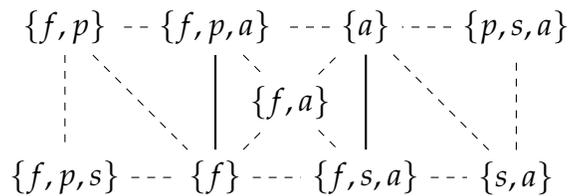
\begin{figure}
\begin{center}
\begin{tikzpicture}
\node (fps) at (-3,0) {$\{f,p,s\}$};
\node (f) at (-1,0) {$\{f\}$};
\node (fsa) at (1,0) {$\{f,s,a\}$};
\node (sa) at (3,0) {$\{s,a\}$};
\node (fp) at (-3,2) {$\{f,p\}$};
\node (fpa) at (-1,2) {$\{f,p,a\}$};
\node (a) at (1,2) {$\{a\}$};
\node (psa) at (3,2) {$\{p,s,a\}$};
\node (fa) at (0,1) {$\{f,a\}$};
\draw[-, dashed] (fps) -- (f) -- (fsa) -- (sa) -- (psa) -- (a) -- (fpa) -- (fp) -- (fps);
\draw[-, dashed] (fp) -- (f) -- (fa) -- (a) -- (sa);
\draw[-, dashed] (fpa) -- (fa) -- (fsa);
\draw[-, thick] (f) -- (fpa);
\draw[-, thick] (a) -- (fsa);
 \end{tikzpicture}
 \end{center}

 \caption{The relative complex from the rooted graph of Figure \ref{fig:graph}}
  \label{fig:graphcomplex}
 \end{figure}

For example, let $G$ be the rooted graph in Figure \ref{fig:graph}, with associated megagreedoid $M$.  The associated relative simplicial complex is given in Figure \ref{fig:graphcomplex}. The triangles are all part of the complex $\Sigma$, every vertex is in $\Gamma$, and every dashed edge is also in $\Gamma$.

 \begin{comment}
 As another example, let $M$ be the megagreedoid from Figure \ref{fig:mega}. Then the tagged intervals in $T$ are the intervals of the form $[x,x]$, intervals $[x,y]$ where $y$ covers $x$, and the intervals $[\{a\}, \{s,c,a\}]$ and $[\{c,a\}, \{s,f,c,a\}]$. The associated relative simplicial complex is given in Figure \ref{fig:complex}. The triangles are simplices, and every vertex and dashed edge is in the subcomplex.
 
\begin{figure}
\begin{center}
\begin{tikzpicture}
\node (fps) at (-3,0) {$\{s,f\}$};
\node (f) at (-1,0) {$\{s\}$};
\node (fsa) at (1,0) {$\{s,c,a\}$};
\node (sa) at (3,0) {$\{c,a\}$};
\node (fp) at (-3,2) {$\{s,f,c\}$};
\node (fpa) at (-1,2) {$\{s,a\}$};
\node (a) at (1,2) {$\{a\}$};
\node (psa) at (3,2) {$\{f,c,a\}$};

\draw[-, dashed] (fp) -- (fps) -- (f) -- (fsa) -- (sa) -- (psa) -- (a) -- (fpa) -- (f) -- (fp);
\draw[-, dashed] (fpa) -- (fsa);

\draw[-, thick] (sa) -- (a) -- (fsa);
 \end{tikzpicture}
 \end{center}

 \caption{The relative complex from the megagreedoid of Figure \ref{fig:mega}}
  \label{fig:complex}
 \end{figure} 
 \end{comment}
\begin{proposition}
Let $M = (A,r)$ be a megagreedoid on $I$, with balanced relative simplicial complex $(\Sigma(A), \Gamma_r(A), \ell)$. Then $F(\Sigma(A), \Gamma_r(A), \ell) = \chi(M) $.
\end{proposition}

 A \emph{shelling} order for $\Gamma \subset \Sigma$ is a linear ordering $F_1, \ldots, F_k$ on the facets of $\Psi$, such that, for each $i$, $\Psi_i \setminus \Psi_{i-1}$ has a unique minimal face, where $\Psi_i = (\cup_{j < i} F_j) \cap \Psi$. 
For shellable, balanced relative simplicial complexes, there is a natural expression for the face quasisymmetric function, in terms of fundamental quasisymmetric functions:
\begin{proposition}
Let $(\Sigma, \Gamma, \rho)$ be a balanced complex, with a shelling order $F_1, \ldots, F_n$. Then $F(\Sigma, \Gamma, \rho) = \sum_{i=1}^n F_{\rho(R(F_i)), n }$, where $R(F_i)$ is the unique minimal face of $\Psi_i \setminus \Psi_{i-1}$.
\end{proposition}

We show that $(\Sigma(A), \Gamma_r(A), \ell)$ is relatively shellable. In the full version, we relate our construction to EC shellability \cite{kozlov}. For the present paper, it is actually shorter to explain the shelling in terms of the greedy algorithm. Fix a linear order on $I$. Consider two $c, c'$ given by $x_1 < \cdots < x_k$ and $y_1 < \cdots < y_{k'}$, let $i$ be the first index where $x_i \neq y_i$. Then we define $c < c'$ if either $r(\{x_1, \ldots, x_{i+1}\}) < r(\{y_1, \ldots, y_{i+1} \})$, or we have equality, and $x_{i+1} < y_{i+1}$. That is, we make a greedy choice. We call this the greedy ordering.
\begin{theorem}
Let $M$ be a megagreedoid. Then the greedy ordering is a shelling order for $(\Sigma(M), \Gamma_r(M))$. Thus $\chi(M)$ is $F$-positive.
\end{theorem}
Note that two types of descents from $A$-feasible permutations come from the greedy ordering: they are cases where there was a `better' local choice for a chain. The other two types of descents come from cases where the facet of $\Sigma$ contains a facet of $\Gamma$.

\section{Hopf monoids and basic quasisymmetric functions}
\label{sec:hopf}

In this section, we dicuss combinatorial Hopf monoids, their characters, and their quasisymmetric functions. Hopf monoids are a generalization of graphs, posets and matroids. The idea is that we have some notion of combinatorial structure, called a species \cite{joyal}. Moreover, we have rules for combining and decomposing these structures in a coherent way. Hopf monoids in species were originally introduced in \cite{aguiar-mahajan-1}. Throughout this section, $A \sqcup B$ denotes disjoint union.

\begin{definition}
A \emph{species} is functor $\spe{F}: \mbox{Set} \to \mbox{Vec}$ from the category of finite sets with bijections, to the category of vector spaces over a field.
For each finite set $I$, $\spe{F}_I$ is a vector space, and for every bijection $\sigma:I \to J$ between finite sets, there is an isomorphism $\spe{F}_{\sigma}: \spe{F}_I \to \spe{F}_J$, such that $\spe{F}_{\sigma \circ \tau} = \spe{F}_{\sigma} \circ \spe{F}_{\tau}$ for every pair $\sigma: I \to J$, $\tau: K \to I$.
 It is \emph{connected} if $\dim \spe{F}_{\emptyset} = 1$. All species in this paper are connected.
\end{definition}

\begin{definition}
A \emph{connected Hopf monoid} is a species $\spe{F}$, equipped with compatible multiplication and comultiplication maps. That is, for every pair of finite sets $S,T$, we have a multiplication map $\mu_{S,T}: \spe{F}_S \otimes \spe{F}_T \to \spe{F}_{S \sqcup T}$ and a comultiplication map $\Delta_{S,T}: \spe{F}_{S \sqcup T} \to \spe{F}_S \otimes \spe{F}_T$. We denote the product of $\spe{f} \in \spe{F}_S$, $\spe{g} \in \spe{F}_T$ by $\spe{f} \cdot \spe{g}$, and we use an analogue of sumless Sweedler notation $\Delta_{S,T}(\spe{f}) = \spe{f}|_S \otimes \spe{f}/S$. Here are some of the axioms:
\begin{enumerate}
\item $(\spe{f} \cdot \spe{g}) \cdot \spe{h} = \spe{f} \cdot (\spe{g} \cdot \spe{h})$ whenever the multiplication is defined.
\item $((\spe{f}|_{B})|_A \otimes (\spe{f}|_{B})/A) \otimes \spe{f} / B = \spe{f}|_A \otimes (\spe{f}/A|_{B - A} \otimes \spe{f}/A/ (B - A))$ for $A \subseteq B \subseteq I$.
\item for $A \sqcup B = S \sqcup T$, and $\spe{f} \in \spe{F}_S, \spe{g} \in \spe{F}_T$, $(\spe{f}\cdot \spe{g})|_A = \spe{f}|_{A \cap S} \cdot \spe{g}|_{A \cap T}$, and $(\spe{f} \cdot \spe{g}) / A = (\spe{f}) / (B \cap S) \cdot (\spe{g}) / (B \cap T)$.
\end{enumerate}
\end{definition}

Given a connected Hopf monoid $\spe{H}$, there is a generalization of M\"obius inversion, and group inversion, called the antipode map. It can be defined recursively: for a finite set $I$, Hopf monoid $\spe{F}$, and $\spe{f} \in \spe{F}_I$, $s_I(\spe{f}) = - \sum_{J \subset I} S_J(\spe{f}|_J) \cdot \spe{f} /J$. This map is linear.

A natural invariant associated to a Hopf monoid is the chromatic quasisymmetric function for a character. 
Given a monoid in species $\spe{M}$, a \emph{character} is a function $\zeta: \spe{M}_I \to \mathbb{K}$ for all $I$, such that for all finite sets $I$, $J$, $\spe{x} \in \spe{M}_I$, $\spe{y} \in \spe{M}_J$, we have $\zeta(\spe{x}) \zeta(\spe{y}) = \zeta(\spe{x}\cdot \spe{y})$. Given a character $\zeta$, and $\spe{x} \in \spe{M}_I$, the basic quasisymmetric function is defined by $$\chi(\spe{x}) = \sum_{\emptyset = S_0 \subset S_1 \subset \cdots \subset S_K = I} (\prod_{i=1}^k \zeta(\spe{x}|_{S_i} / S_{i-1})) M_{\{|S_1|, |S_2|, \cdots, |S_{k-1}|\}, |I|}.$$

Given a Hopf monoids $\spe{H}$, a Hopf submonoid of $\spe{H}$ is a species $\spe{K}$ where $\spe{K}_I \subseteq \spe{H}_I$, and is preserved under multiplication and comultiplication.

\begin{example}
Here we detail the Hopf monoid of megagreedoids $\spe{MG}$.
We let $\spe{MG}_I$ be the vector space with basis given by megagreedoids. Given sets $I = S \sqcup T$, we define $\mu_{S,T}: \spe{MG}_S \otimes \spe{MG}_T \to \spe{MG}_{S \sqcup T}$ 
by sending a pair of megagreedoids to their direct sum, and extending by linearity. 
Likewise, we define $\Delta_{S,T}(M) = M|_S \otimes M / S$ if $S$ is feasible, and $0$ otherwise. Extend by linearity to get a map $\Delta_{S,T}: \spe{MM}_{S \sqcup T} \to \spe{MM}_S \otimes \spe{MM}_T$. This turns $\spe{MG}$ into a Hopf monoid. 
In the full version of the paper, we show that $\chi(A,r)$ is the quasisymmetric function that arises from the basic character $\zeta$ given by $\zeta(A,r) = 1$ if $A$ is boolean and $r$ is modular, and $0$ otherwise.

\end{example}
We do not have a simple, cancellation free formula for the antipode. Is there one? We only can manage to use Takeuchi's formula to get our combinatorial reciprocity result, but it would be nice to have the entire antipode.

\begin{example}
For the polymatroid species $\spe{PM}$, the basis of $\spe{PM}_I$ consists of all polymatroids on $I$. This is a Hopf submonoid of $\spe{MG}$, and is isomorphic to the Hopf monoid of generalized permutahedra \cite{aguiar-ardila}. It is related to several other Hopf monoids, coming from hypergraphs, building sets, and simplicial complexes.
\end{example}

\begin{comment}
\begin{example}
Our next example is the Hopf monoid of posets with inconsistent pairs $\spe{PIP}$. The basis of $\spe{PIP}_I$ consists of all triples $(\spe{p}, \spe{I}, w)$ where $\spe{p}$ is a poset with inconsistent pairs $\spe{I}$, and $w: I \to \mathbb{N}$ is a weight function. The multiplication on $\spe{P}$ is induced by disjoint union of partial orders, and disjoint union of the inconsistent pairs.
Given $\spe{p}$, and $S \subseteq I$, we let $\spe{p}[S]$ denote the induced subposet. Then $\spe{p}|_S = \spe{p}[S]$ provided $\spe{p}[S]$ is an order ideal of $\spe{p}$. Similarly, if $\spe{p}[S]$ is a lower order ideal, we let $\spe{p}/S = \spe{p}[I \setminus S]$. We define $\spe{I}|_{S}$ (respectively, $\spe{I}/S$) to consist of those inconsistent pairs where both endpoints are in $S$ (respectively, in $T$). We define $w|_S$ bey $w|_S(v) = w(v) + $ the number of inconsistent pairs $(v,w)$ where $w \in I \setminus S$. Similarly, $w/S$ is the restriction of $w$ to $I \setminus S$. This turns $\spe{PIP}$ into a Hopf monoid.

The basic character on $\spe{P}$ is given by $\zeta(\spe{p}) = 1$ if $\spe{p}$ is an antichain with no inconsistent pairs, and $0$ otherwise. The resulting quasisymmetric function $\chi(\spe{p})$ is a generalization of Gessel's $P$-partition quasisymmetric function which counts strict order-preserving maps $f: \spe{p} \to \mathbb{N}$, where each map is weighted by $\prod_{i \in I} x_{f(i)}$. 
The map sending a poset $\spe{p}$ to its megagreedoid is also a morphism of Hopf monoids.
\end{example}
\end{comment}

\begin{example}
Our next example is the Hopf monoid of posets $\spe{P}$. The basis of $\spe{P}_I$ consists of all partial orders on $I$. The multiplication on $\spe{P}$ is induced by disjoint union of partial orders.
Given $\spe{p}$, and $S \subseteq I$, we let $\spe{p}[S]$ denote the induced subposet. Then $\spe{p}|_S = \spe{p}[S]$ provided $\spe{p}[S]$ is an order ideal of $\spe{p}$. Similarly, if $\spe{p}[S]$ is a lower order ideal, we let $\spe{p}/S = \spe{p}[I \setminus S]$.  When $S$ is not a lower order ideal, then $\spe{p}|_S = \spe{p}/S = 0$. This turns $\spe{P}$ into a Hopf monoid.

The basic character on $\spe{P}$ is given by $\zeta(\spe{p}) = 1$ if $\spe{p}$ is an antichain, and $0$ otherwise. The resulting quasisymmetric function $\chi(\spe{p})$ is Gessel's $P$-partition quasisymmetric function which counts strict order-preserving maps $f: \spe{p} \to \mathbb{N}$, where each map is weighted by $\prod_{i \in I} x_{f(i)}$. 
The map sending a poset $\spe{p}$ to its megagreedoid embeds $\spe{P}$ as a Hopf submonoid of $\spe{MG}$. Thus, we obtain $F$-positivity from our more general result.
\end{example}

\begin{example}
Another Hopf submonoid is given by the greedoid species $\spe{Gr}$, where the basis of $\spe{Gr}_I$ consists of all greedoids on $I$. 
Recall that contraction for greedoids is only defined for \emph{rank-feasible} sets. Given a rank feasible subset $R$ for a greedoid $\spe{g}$ with rank function $r$, the restriction $\spe{g}|_R$ has rank function $r|_R: 2^R \to \mathbb{N}$ given by $r|_R(S) = r(S)$ for all $S \subseteq R$. The contraction $\spe{g}/R$ has rank function $r/R:2^{I \setminus R} \to \mathbb{N}$, where $r/R(S) = r(S \cup R) - r(R)$. Much like for megagreedoids, restriction and contraction can be used to define a comultiplication: $\Delta_{S,T}(\spe{g}) = \spe{g}|_S \otimes \spe{g}/S$, where the tensor is $0$ if $S$ is not rank-feasible. Extend by linearity to get a map $\Delta_{S,T}: \spe{Gr}_I \to \spe{Gr}_S \otimes \spe{Gr}_T$. Using direct sum of greedoids as a multiplication operation, we see that $\spe{Gr}$ is a Hopf monoid. Moreover, our construction $\spe{g} \mapsto \spe{m}_{\spe{g}}$ associating a megagreedoid to greedoid embeds $\spe{Gr}$ as a Hopf submonoid of $\spe{MG}$. 
There are lots of Hopf submonoids of greedoids: matroids, posets, antimatroids, interval antimatroids, and directed branching greedoids.

\end{example}

\begin{figure}[t]

\begin{center}
\begin{tikzpicture}
\node (e) at (0,0) {$c$};
\node (f) at (-1,1) {$f$};
\node (n) at (1,1) {$a$};
\node (u) at (-1,2) {$p$};
\node (s) at (1,2) {$s$};

\draw[-] (f) -- (e);
\draw[-] (n) -- (e);
\draw[-] (u) -- (f);
\draw[-] (f) -- (n);
\draw[-] (u) -- (s);
\draw[-] (s) -- (n);

\end{tikzpicture}
\begin{tikzpicture}
\node (e) at (0,0) {$c$};
\node (f) at (-1,1) {$f$};
\node (n) at (1,1) {};
\node (u) at (-1,2) {$p$};
\node (s) at (1,2) {};

\draw[-] (f) -- (e);

\draw[-] (u) -- (f);
\draw[-] (f) -- (n);
\draw[-] (u) -- (s);

\end{tikzpicture}
\begin{tikzpicture}
\node (e) at (0,0) {$c$};

\node (n) at (1,1) {$a$};

\node (s) at (1,2) {$s$};

\draw[bend left, -] (n) to (e);
\draw[-] (n) to (e);
\draw[-] (e) -- (s);

\draw[-] (s) -- (n);

\end{tikzpicture}

\end{center}
\caption{A rooted graph, a restriction, and a contraction.}
\label{fig:ex}
\end{figure}
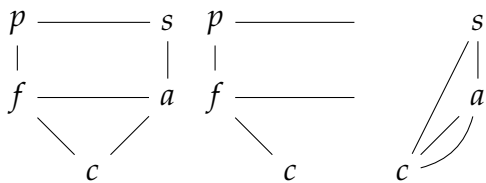

\begin{example}
Finally, we discuss the Hopf monoid of rooted connected graphs $\spe{RG}$, and the chromatic symmetric function of a rooted connected graph. We show that $\spe{RG}$ is a Hopf submonoid of $\spe{MM}$.

The basis of $\spe{RG}_I$ is given by all rooted connected graphs $\spe{g}$ with vertex set $I \sqcup \{r \}$. Note that we do allow multiple edges. Given two rooted graphs $\spe{g}$ and $\spe{h}$ such that $\spe{g} \cap \spe{h} = \{r\}$,  let $\spe{g} \cdot \spe{h}$ be their union: note that $r$ is the root of all three graphs. 

Next, we describe the comultiplication. We only define the restriction and contraction for subsets $S \cup \{r\} \subset I \cup \{r\}$ that induce connected subgraphs. The restriction map $\spe{g}|_S$ is the subgraph on $S \cup \{r\}$ which contains any edge that has at least one endpoint in $S$. Note that some edges become half-edges in the restriction.
On the other hand, $\spe{g} / S$ is given by contracting $S \cup \{r\}$ into one vertex, labeled $r$, and regarded as the root of $\spe{g} / S$. Please see Figure \ref{fig:ex} for an example. Here we restrict and contract with respect to $\{f,p \}$. If $S \cup \{r\}$ is not connected, we define $\spe{g}|_S = \spe{g}/S = 0$. This defines a comultiplication $\Delta_{S,T}: \spe{RG}_{S \sqcup T} \to \spe{RG}_S \otimes \spe{RG}_T$
given by $\Delta_{S,T}(\spe{g}) = \spe{g}|_S \otimes \spe{g}/S$.

The basic character on $\spe{G}$ is given by $\zeta(\spe{g}) = 1$ if $\spe{g}$ is a star graph, centered at $r$, and $0$ otherwise. Then $\chi(\spe{g})$ is our chromatic quasisymmetric function.

There is a natural map $\spe{g} \mapsto \spe{m}_{\spe{g}}$ sending a rooted graph to its megagreedoid.
This embeds $\spe{RG}$ as a Hopf submonoid of $\spe{MG}$.
There is also an interesting Hopf submonoid consisting of those graphs for which $r$ is adjacent to all vertices. This is isomorphic to the Hopf monoid of graphs appearing in Aguiar and Ardila \cite{aguiar-ardila}.

\end{example}

\section{Conclusion}

We end with lots of open questions. First, is there a notion of duality for megagreedoids? Also, is there an analogue of the $G$-invariant introduced by Dersken and Fink \cite{derksen-fink}? Is there a notion of valuative invariant? What about a Tutte polynomial? Are there other interesting examples beyond rooted connected graphs?

Regarding the quasisymmetric functions themselves, there is another obvious invariant, where we sum over all $r$-feasible functions that are minimized at a unique basis. Is the resulting invariant $F$-positive? It appears to be so, but our labeling does not prove shellability in this case. Also, what can we say about the chromatic quasisymmetric functions of Shareshian and Wachs?

Given a megagreedoid, with relative complex $(\Gamma, \Delta)$, is $\Gamma$ shellable? Does it have a convex ear decomposition?
What more can be said about coloring rooted connected graphs? Are there nice formulas for special graphs?

%\small
\bibliographystyle{abbrv}  
%\bibliographystyle{amsalpha}  
%\nocite*
%\bibliographystyle{abbrvnat}
%\bibliographystyle{amsplain}  
\bibliography{hopf_monoid_coloring}

\end{document}